\documentclass[12pt]{amsart}
\usepackage{amsfonts, amsmath, amsthm}

\newtheorem{pro}{Proposition}[section]
\newtheorem{thm}[pro]{Theorem}
\newtheorem{lem}[pro]{Lemma}

\newtheorem{ass}[pro]{Assertion}

\newtheorem{rmkks}[pro]{Remarks}

\newtheorem{cor}[pro]{Corollary}

\newtheorem*{qthm}{Theorem}

\theoremstyle{definition}
\newtheorem{dfn}[pro]{Definition}
\newtheorem{dfns}[pro]{Definitions}

\theoremstyle{remark}
\newtheorem*{rmk}{Remark}
\newtheorem*{rmks}{Remarks}

\newcommand{\pmm}{primitive meridian}

\newcommand{\scc}{simple closed curve}
\newcommand{\s}{\Sigma}

\newcommand{\sft}{{swallow follow torus}}

\newcommand{\eg}{{\it e.g.}}

\newcommand{\del}{\partial}

\newcommand{\hhs}{Heegaard surface}

\title{Morimoto's Conjecture for m-small knots}

\date{\today}
\address{Department of Mathematics, Nara Women's University
Kitauoya Nishimachi, Nara 630-8506, Japan}
\address{Department of mathematical Sciences, University of
Arkansas, Fayetteville, AR 72701}
\email{tsuyoshi@cc.nara-wu.ac.jp} \email{yoav@uark.edu}
\author{Tsuyoshi Kobayashi}
\author{Yo'av Rieck}
\thanks{Both authors are supported in part by JSPS grants;  the
second named author was a JSPS fellow P00024.}

\begin{document}

%\subjclass{}%
%\keywords{}%

%\date{}%
%\dedicatory{}%
%\commby{}%
% ----------------------------------------------------------------
\begin{abstract}

Let $X$ be the exterior of connected sum of knots and $X_i$ the
exteriors of the individual knots.  In \cite{morimoto1} Morimoto
conjectured (originally for $n=2$) that $g(X) < \s_{i=1}^n g(X_i)$
if and only if there exists a so-called \em primitive meridian \em
in the exterior of the connected sum of a proper subset of the
knots. For m-small knots we prove this conjecture and bound the
possible degeneration of the Heegaard genus (this bound was
previously achieved by Morimoto under a weak assumption
\cite{morimoto2}):

$$\s_{i=1}^n g(X_i) - (n-1) \leq g(X) \leq \s_{i=1}^n g(X_i).$$

\end{abstract}
\maketitle
% ----------------------------------------------------------------
\section{Introduction}

This proceeding article is based on a talk given in Waseda
University on the 18 December 2002, about Morimoto's Conjecture
which is concerned with the behavior of Heegaard genus of knot
exteriors under connected sum. For a knot exterior, we consider
two equivalent decomposition: the first is given by \em tunnel
system, \em which is a collection of embedded arcs and the other a
\em Heegaard splitting, \em given by an embedded surface. (See
next section for standard definition.) The complexity of a tunnel
system is the number of arcs and the complexity of a Heegaard
surface is its genus.  The complexity of a knot $K \subset M$ is
the minimal number of tunnels required for a tunnel system (called
the tunnel number and denoted $t(K)$) or the genus of the minimal
genus Heegaard surface (called the Heegaard genus and denoted
$g(X)$, here $X = M \setminus N(K)$). It is immediate from
definitions that $g(X) = t(K) + 1$. Let $K_1 \subset M_1$ and $K_2
\subset M_2$ be two knots, and let $K (=K_1 \# K_2) \subset M (=
M_1 \# M_2)$ be the connected sum. In Section \ref{sec:background}
we recall the easy fact: $t(K) \leq t(K_1) + t(K_2) + 1$.  (More
generally, $t(K) \leq \Sigma _{i=1}^n t(K_i) + (n-1).$)
Translating this into the language of Heegaard genus, we get $g(X)
\leq g(x_1) + g(X_2)$. (More generally, $g(X) \leq \Sigma _{i=1}^n
g(X_i).$)  In this paper we use the notation of Heegaard genus
which seems simpler than tunnel number. Note that a knot exterior
$X$ (resp. $X_i$) has a distinguished slope, the meridian, denoted
$\mu$ (resp. $\mu_i$). Since we are interested in the knot
exterior, instead of studying the connected sum of knots, we
consider what happens to the exteriors. It is easy to see that $X
= X_1 \cup_A X_2$, where $A \subset \del X_i$ is an annulus that
is a neighborhood of $\mu_i$ in $\del X_i$ ($i=1,2$). We denote
that operation $X_1 \star X_2$ but emphasize that $X_1 \star X_2$
depends not only on $X_i$ but on $\mu_i$ as well.  It should be
noted that the operation $\star$ is more closely related to torus
decomposition than connected sum.

Y. Moriah and J.H. Rubinstein \cite{moriah-rubinstein} showed that
there exist knots $K_1$ and $K_2$ for which $g(X) = g(X_1) +
g(X_2)$.  K. Morimoto \cite{morimoto1} constructed examples for
which $g(X) = g(X_1) + g(X_2) - 2$ and T. Kobayashi
\cite{kobayashi-degeneration} generalized them to examples were
$g(X) = g(x_1) + g(X_2) - n$ for arbitrarily large $n$.  These and
all other known examples have one thing in common: at least one of
the two knot exteriors has a minimal genus Heegaard surface with a
\pmm.

\begin{dfn}
\label{dfn:pm}

A minimal genus Heegaard splitting for a knot exterior $X$ has a
{\it\pmm} if there exists a compressing disk for the compression
body that intersect a spanning annulus with a meridian slope
exactly once.

\end{dfn}

\begin{rmk}
\label{rmk:pm-is-min-genus}

A non-minimal genus Heegaard surface, even if it is irreducible,
is never said to have a \pmm.  Otherwise we need to modify
Assertion \ref{sufficient}.

\end{rmk}

Morimoto further observes that the existence of a \pmm\ guarantees
a degeneration of the genus:

\begin{ass}
\label{sufficient}

Let $X_1$ and $X_2$ be two knot exteriors and let $X = X_1 \star
X_2$ be the exterior of the connect sum.

If $X_1$ or $X_2$ contains a \pmm\ then $g(X) \leq g(x_1) + g(X_2)
- 1.$

\end{ass}

\begin{proof}

Recall that $X_1 \star X_2$ is obtained from $X_1$, $X_2$ by
identifying meridional annuli.

Let $\s_1 \subset X_1$ and $\s_2 \subset X_2$ be minimal genus
Heegaard surfaces, say $\s_2$ contains a \pmm. After gluing $X_1$
to $X_2$ we surger the Heegaard surfaces together along an annulus
that runs across $A$. The component of $X_1 \star X_2$ containing
the boundary is easily seen to be a compression body: it is simply
a neighborhood of the boundary union neighborhoods of tunnels for
both knots.  The other component is obtained by gluing the
handlebody components of $X_i$ cut open along $\s_i$ to each other
along an annulus, which is general does not yield a handlebody.
However, since $X_2$ contained a \pmm\ the annulus is longitudinal
there, and a handlebody is obtained.
\end{proof}

Before stating our results we define:

\begin{dfns}
\label{dfn:msmall}
\begin{enumerate}

\item A surface with non empty boundary, properly embedded in a
manifold is called \em essential \em if it is incompressible,
boundary incompressible and not boundary parallel.

\item A knot exterior is called \em meridionally small \em if
there is no essential surface with non empty boundary whose
boundary forms parallel copies of the meridian; in other words the
meridian is not a boundary slope.

\end{enumerate}
\end{dfns}

We prove two results, the first numerical:

\begin{thm}[The Numerical Theorem]
\label{thm:numerical}

Let $\{X_i\}_{i=1}^n$ be exteriors of m-small knots and $X =
\star_{i=1}^n X_i$ be the exterior of the connected sum. Then:

$$\s_{i=1}^n g(X_i) - (n-1) \leq g(X) \leq \s_{i=1}^n g(X_i).$$

\end{thm}

\begin{rmks}{\rm
The assumption required for this bound is in  fact weaker than
m-smallness: we just need a minimal genus Heegaard surface to
weakly reduce to a swallow follow torus. (see Theorems
\ref{thm:basic} and \ref{thm:bound-on-degeneration}). }
\end{rmks}

Thus, in contrast to the examples of Morimoto and of Kobayashi
mentioned above, the degeneration of Heegaard genus under
connected sum is bounded.  This result, with a small assumption
(that none of the ambient manifolds has a lens space component)
was obtained by Morimoto \cite{morimoto2}.

In this theorem and throughout the paper we are considering knot
exteriors.  The statements can be rephrased in the language of
manifolds with boundary torus. Let $(X_i,\mu_i)$ be the manifold
$X_i$ with a choice of meridian $\mu_i$ ($i=1,..,n)$.  By
Hatcher's Theorem \cite{hatcher} $(X_i,\mu_i)$ is $\mu_i$-small
for all but finitely many choices of $\mu_i$.  Therefore for all
other values of $\mu_i$ the bound above is valid; we get the
following theorem (other theorems in this paper can be modified
similarly).

\begin{qthm}[The Numerical Theorem (\ref{thm:numerical}) version 2]
Given $X_i$ ($i=1,..,n$) manifolds with boundary torus, after
excluding a finite set of slopes from each $X_i$ we get the
following inequality for all remaining slopes:

    $$\s_{i=1}^n g(X_i) - (n-1) \leq \star_{i=1}^n (X_i,\mu_i))
    \leq \s_{i=1}^n g(X_i).$$
\end{qthm}

Our second result is geometric; it generalizes Morimoto's
Conjecture for the connected sum of $n$ knots (we return to the
language of knot exteriors):

\begin{thm}[Morimoto's Conjecture]
\label{thm:geometric}

Let $\{X_i\}_{i=1}^n$ be exteriors of m-small knots. Then:

If $g(X) < \s_{i=1}^n g(X_i)$ there exists some $I \subset
\{1,..,n\}$ a proper subset so that $X = \star_{i \in I} X_i$
contains a \pmm.
\end{thm}

We remark that by saying ``$I$ is a proper subset'' we mean that
$I \neq \{1,..,n\}$ and $I \neq \emptyset$.

This theorem too has a second version, again using Hatcher's
Theorem:

\begin{qthm}[Morimoto's Conjecture (\ref{thm:geometric}) version 2]
Given $X_i$ ($i=1,..,n$) manifolds with boundary torus, after
excluding a finite set of slopes from each $X_i$ for all remaining
slopes we get:

If $ g(\star_{i=1}^n (X_i,\mu_i))  < \s_{i=1}^n g(X_i)$ there
exists some $I \subset \{1,..,n\}$ a proper subset so that $X =
\star_{i \in I} X_i$ contains a \pmm.
\end{qthm}

The main tool for proving both these results is:

\begin{thm}[The Swallow Follow Theorem]
\label{thm:basic}

Let $X_i$ be a collection of exteriors of meridionally small
knots.

Then any Heegaard splitting of $X$ (the exterior of the connected
sum) weakly reduces to a swallow follow torus.

\end{thm}

{\bf Acknowledgement:} We would like to thank Kanji Morimoto for
helpful conversations.  The second named author: this research was
conducted while I was a JSPS fellow at Nara Women's University.  I
would like to thank the university, the math department and
particularly Tsuyoshi Kobayashi for wonderfully warm hospitality,
and Kouki Taniyama for giving me the opportunity to visit Waseda
University.

\section{Background}
\label{sec:background}

We review some standard definitions (see also \cite{hempel} or
\cite{jaco}): let $M$ be a closed 3-manifold. A knot $k \subset M$
is a smooth embedding of $S^1$ into $M$.  The knot exterior $X$
(which is a manifold with boundary torus) is $M \setminus N(k)$,
where $N(\cdot)$ is an open normal neighborhood.  A \em
compression body \em is a manifold $C$ with distinguished boundary
component (denoted $\del_+ C$, and $\del_- C = \del C \setminus
\del_+ C$) so that after compressing $\del_+$ maximally a
collection of balls and components of the form ($\del_- C \times
[0,1]$) are obtained. A compressing disk for $\del_+ C$ is called
a meridian disk.  A \em handlebody \em (denoted $H$) is the
special case of compression body where only balls are obtained
(equivalently if $\del_- H=\emptyset$).  A \em tunnel system \em
is a collection of arcs properly embedded in the knot exterior
whose exterior is a handlebody.  It is well known (and an easy
consequence of Morse theory) that every knot has a tunnel system.
The tunnel number $t(K)$ is the least number of arcs required for
a tunnel system. If $T$ is a tunnel system, $\del(N(\del X \cup
T))$ is a surface that decomposes the knot exterior into a
compression body $C$ (with $\del_- C = \del X$) and a handlebody
$H$. A closed surface that decomposes $X$ into a handlebody and a
compression body is called a \em Heegaard surface \em (denoted
$\s(X)=\s$).  It is well known that any Heegaard surface is given
as the neighborhood of some tunnel system.  Note that the genus of
the $\s$ (denoted $g(\s)$) is exactly one more than the number of
tunnels.

A Heegaard splitting is called stabilized if there are meridional
disks on opposite sides intersecting exactly once, non-stabilized
otherwise.  If a Heegaard splitting is stabilized, by cutting
along one of the disks a Heegaard splitting of lower genus is
obtained (it is a consequence of Scharlemann and Thompson
\cite{scharl-thomp} that the converse is also true).  We will be
mostly concerned with minimal genus Heegaard splittings which are
therefore always non-stabilized.

An extremely important concept for the study of Heegaard
splittings is \em strong irreducibility \em (see
\cite{casson-gordon}).  A \em weak reduction \em for $\s$ is a
pair meridian disks on opposite sides of $\s$ that are disjoint.
$\s$ is called strongly irreducible if there does not exist a weak
reduction, weakly reducible otherwise.  If after compressing $\s$
to both sides an essential surface $F$ is obtained we say that $F$
was obtained from $\s$ by weak reduction.  Casson and Gordon
\cite{casson-gordon} showed that if a non stabilized Heegaard
surface weakly reduces, \em some \em weak reduction yields an
essential surface.  This was the first time a Heegaard surface was
used to produce an essential surface, but by no means the last.

A \em connected sum \em of knots $K_i \subset M_i\,(i=1,2)$ is
defined (much like the case $M_i \cong S^3$) by removing a small
ball around each knot and gluing the spheres obtained so that the
endpoint of the arcs match up. After removing the intersection of
$N(K)$ with the sphere an essential annulus is obtained, called a
\em decomposing annulus. \em  The connected sum is denoted $K_1 \#
K_2$ and is naturally a knot in $M_1 \# M_2$. It is now straight
forward to see that $t(K) \leq t(K_1) + t(K_2) + 1$ by taking the
union of tunnel system for $K_1$ and $K_2$ and a tunnel that is an
essential arc on the decomposing annulus. As noted above the
exterior of $K_1 \# K_2$ is obtained from the exteriors of $K_1$
and $K_2$ by gluing a meridional annulus.  We denote the exterior
of $K_1 \# K_2$ by $X_1 \star X_2$ (or just $X$), where $X_i$ is
the exterior of $K_i$. Connected sum of $n$ knots is defined by
induction, and similar notation is used.

Scharlemann and Thompson \cite{scharl-abby} refined Casson and
Gordon's weak reduction \cite{casson-gordon} in a construction
called \em untelescoping \em (Kobayashi \cite{kobayashi} showed
that untelescoping is strictly finer than weak reduction).  We
briefly describe untelescoping: by Casson and Gordon any
non-stabilized Heegaard splitting can be untelescoped to give a
collection of closed, essential surfaces (if the Heegaard surface
was strongly irreducible this collection is empty and we are
done). Each submanifold obtained by cutting the manifold open
along these surfaces inherit a Heegaard splitting (called the \em
induced \em Heegaard splitting).  If the induced Heegaard
splitting weakly reduces, apply Casson and Gordon's result again
(see Sedgwick \cite{eric-agt} for the case with boundary) to
obtain a larger collection of essential surfaces. We repeat this
process that eventually terminates. Scharlemann and Thompson show
that the induced Heegaard splitting after this final stage is
strongly irreducible.  We will be using these facts in an
essential way in the proof presented in Section \ref{sec:morimoto}
refer the reader to \cite{scharl-abby} and \cite{kobayashi} for a
more detailed account of untelescoping. Since getting $\s$ back
from the induced Heegaard splittings is the converse of
untelescoping, called \em amalgamating. \em

\section{basic constructions}
\label{sec:basic}

The construction behind all our work is based on the following
theorem:

\begin{qthm}[\ref{thm:basic}]

Let $X_i$ be a collection of exteriors of meridionally small
knots.

Then any Heegaard splitting of $X$ (the exterior of the connected
sum) weakly reduces to a swallow follow torus.

\end{qthm}

\begin{proof}[Sketch of proof]

Let $\s$ be a Heegaard surface for $X$.  Let $A =
\{A_j\}_{j=1}^{n-1}$ be some collection of essential annuli so
that $X$ cut open along $A$ is $\bigsqcup_{i=1}^n X_i$
($\bigsqcup$ denotes disjoint union).  We consider three cases.

\begin{description}

\item[Case One: assume $\s$ is strongly irreducible] We show that
this leads to a contradiction.  If $\s$ is strongly irreducible we
can then isotope it to intersect $A$ in curves that are all
essential in both $A$ and $\s$ (this is a standard application of
strong irreducibility and we omit the details). For $i=1,2$ let
$\s_i$ be the collection of surfaces obtained by cutting $\s$
along $A$, with $\s_i \subset X_i$ ($\s_i$ need not be connected).
Since $\s$ is a Heegaard surface for $X$ (and $X$ is not a
compression body) $\s$ compresses into both sides.  Let $D_+$ and
$D_-$ be compressing disks for $\s$ on opposite sides of $\s$. By
minimizing the intersection of $D_+$ and $D_-$ with $A$, and using
innermost and outermost disk arguments, we obtain a compression or
a boundary compression from some component of $\s_i$, and since
$\partial X_i$ is a torus, boundary compression implies
compression (this point requires a little care about components
that are boundary parallel annuli, but we ignore that here).  We
choose a normal direction for $\s$.  Each component of $\s_1$ and
$\s_2$ inherits a normal direction.  Some component of some $\s_i$
compresses to one side and some component of some $\s_j$ (possibly
$i=j$) compresses to the other side. A compression for $\s_i$ is
also a compression for $\s$ (recall that all curves of
intersection are essential in both). By strong irreducibility the
compressing disks must intersect, and therefore both compress the
same component of the same $\s_i$. Furthermore strong
irreducibility implies that all other components of $\s_1$ and
$\s_2$ are essential. This violates m-smallness of the knot
exteriors $X_i$.  (Compare this with \cite{rieck-sedgwick1} and
\cite{rieck-sedgwick2}: there too the authors produce bounded
essential surfaces using Heegaard surfaces.)

\end{description}

Thus we may assume that $\s$ weakly reduces, and we can maximally
untelescope it (\cite{scharl-abby}, recall Section
\ref{sec:background}) obtaining a collection of closed essential
surfaces denoted $S$, with $S_j$ denoting a component of $S$.
Minimize the intersection $S$ with $A$. Any component of $S$ that
actually intersects $A$ is broken up into pieces that are
incompressible, and if some such component is also of negative
Euler characteristic, it must be essential, a contradiction
similarly to above. (Note that no such component is a disk.) But
this would contradict m-smallness. We conclude that each component
$S_j$ of $S$ is exactly one of the following:

\begin{enumerate}
\item $S_j \cap A = \emptyset$, or---

\item each component of $S_j$ cut open along $A$ is an annulus
(these annuli must be boundary parallel).
\end{enumerate}

We are now ready for the remaining two cases:

\begin{description}

\item[Case Two] $(\forall j) S_j \cap A = \emptyset$.  Cutting $X$
open along $S$ we get a collection of manifolds, one of them
containing $A$. By the maximality of the untelescoping the
Heegaard splitting for this component is strongly irreducible, but
in Case One of this proof we showed that if an annular manifold
contains a strongly irreducible Heegaard surface then it contains
an essential surface with boundary on $A$. (In our main paper we
show that that proof is valid here as well.) This surface remains
essential in $X$, and can be used to contradict m-smallness

\item[Case Three] $(\exists j) S_j \cap A \neq \emptyset$. We
following an argument of Morimoto \cite{morimoto2}.  $S_j$ is cut
up to annuli, since any component of negative Euler characteristic
would contradict m-smallness. Thus $S_j$ is a torus. To conclude
this sketch we show that a torus that cannot be disjoint from $A$
is a \sft. After arranging the intersection of $S_j$ and $A$ to be
essential and minimal, an outermost subannulus of $A$ cut open
along $S_j$ is an annulus with one boundary on $S_j$ and the other
a meridian. Surgering the torus $S_j$ along this subannulus we
obtain a meridional annulus that must be decomposing. Thus $S_j$
was obtained from some decomposing annulus by tubing along the
boundary and is therefore a swallow follow torus by definition.
This completes the proof.

\end{description}
\end{proof}

We obtain the following corollary, see \cite{morimoto2} for a
detailed description.  This corollary is a restatement of Theorem
\ref{thm:basic} in an easy to use manner.

\begin{dfn}
\label{dfn:hat} Let $Y$ be a knot exterior.  Then $\widehat Y$
denotes the manifold obtained from $Y$ by drilling out a curve
parallel to a meridian on $\del Y$.
\end{dfn}

\begin{cor}
\label{cor:swallow-follow}

Let $T$ be the swallow follow torus obtained in Theorem
\ref{thm:basic}.

Suppose $T$ follows $X_I = \star_{i \in I}X_i$ for some $I \subset
\{1,..,n\}$ and swallows $X_J = \star_{i \in J}X_i$ (here $J =
\{1,..n\} \setminus I$ and of course $I \neq \emptyset$ and $J
\neq \emptyset$).

Then $X \cong X_I \cup_T \widehat X_J$.
\end{cor}

Note in the corollary above $\widehat X_J$ has two boundary
components and $X_I$ only one.  As we shall see, the \pmm\ we are
trying to find is in $\widehat X_J$ and not necessarily in $X_I$,
so \em there is no symmetry between the roles of $X_I$ and $X_J$.
\em More about that later.

\section{numerical bounds}

Corollary \ref{cor:swallow-follow} holds whenever a minimal genus
Heegaard surface weakly reduces to a \sft.  We demonstrate its
usefulness.  In the following theorem we are not assuming that
$X_1$ and $X_2$ are m-small (or even prime).

\begin{thm}
\label{thm:bound-on-degeneration}

Suppose a minimal genus Heegaard surface of $X_1 \star X_2$ weakly
reduces to a \sft\ that follows one and swallows the other.

Then $g(X_1) + g(X_2) - 1 \leq g(X_1 \star X_2)\leq g(X_1) +
g(X_2).$
\end{thm}

\begin{rmk}
Recall that the right hand side inequality is easy and always
holds.
\end{rmk}

Before sketching the proof of Theorem
\ref{thm:bound-on-degeneration} we bring two corollaries. The
first is obtained as follows: let $X$ be the exterior of a
connected sum of m-small knots.  Then $X_I$ and $X_J$ described in
Corollary \ref{cor:swallow-follow} are exteriors of connected sum
of m-small knots as well (perhaps just one summand). By induction
we get:

\begin{qthm}[\ref{thm:numerical}]
\label{cor:numerical}

Let $\{X_i\}_{i=1}^n$ be exteriors of m-small knots and $X =
\star_{i=1}^n X_i$ be the exterior of the connected sum. Then:

$$\s_{i=1}^n g(X_i) - (n-1) \leq g(X) \leq \s_{i=1}^n g(X_i).$$

\end{qthm}

Another corollary is about Morimoto's examples \cite{morimoto2}.
Morimoto has examples of knot exteriors so that $g(X_1) = 2$,
$g(X_2) = 3$ and $g(X_1 \star X_2)$ = 3.  We get:

\begin{cor}
\label{cor:morimoto-example} Let $K_1$ and $K_2$ be knots with
irreducible and a-toroidal exteriors.

If  $g(X_1) + g(X_2) \geq 5$ and $g(X_1 \star X_2) = 3$ then any
minimal genus Heegaard surface for $X_1 \star X_2$ is strongly
irreducible.
\end{cor}

The knots $K_1$ and $K_2$ in Morimoto's example are hyperbolic
knots in $S^3$ and therefore fulfil the assumptions of the
corollary above, thus we get an example of a connected sum that
has only strongly irreducible minimal genus Heegaard surfaces
(this was obtained independently by Y. Moriah
\cite{moriah-example} using different techniques).

\begin{proof}[Proof of Corollary \ref{cor:morimoto-example}]
Assume for contradiction that a minimal genus Heegaard surface
weakly reduces.   Then by Casson and Gordon \cite{casson-gordon}
some weak reduction yields an essential surface $F$.  The genus of
$F$ is at most $g(X_1 \star X_2) - 2$, hence at most 1.  Our
assumptions imply that $X_1 \star X_2$ contains no essential
sphere, hence $F$ must be a torus. If $F$ can be isotoped to be
disjoint from the decomposing annulus it must be parallel to $\del
X_1$ or $\del X_2$, thus a swallow follow torus. On the other hand
if $F$ cannot be made disjoint from the decomposing annulus, we
can arrange the intersection of $F$ and the annulus to consist of
a non-empty collection of essential curves. As we saw above,
surgering $F$ along an outermost subannulus yields an essential
annulus and therefore $F$ is a swallow follow torus. In both cases
we conclude that the minimal genus Heegaard surface weakly reduces
to a \sft.

Hence by Theorem \ref{thm:bound-on-degeneration} the genus reduces
by one at most, contradiction.
\end{proof}

\begin{proof}[Sketch of proof of Theorem \ref{thm:bound-on-degeneration}]

Corollary \ref{cor:swallow-follow} gives the decomposition $X =
X_1 \cup \widehat X_2$ or $X = \widehat X_1 \cup X_2$ (say the
former). Recall that in Theorem \ref{thm:bound-on-degeneration} we
assumed the \hhs\ for $X$ is minimal genus.   Consider the induced
\hhs s on $X_1$ and $\widehat X_2$. If the induced \hhs s on $X_1$
and $\widehat X_2$ are not minimal genus, replacing them with a
lower genus \hhs\ and amalgamating we get a lower genus \hhs\ for
$X$, contradiction. Using these \hhs s we get (\eg\
\cite{schultens}):

\begin{equation}\label{equ:amalgamation}
g(X) = g(X_1) + g(\widehat X_2) - g(T) = g(X_1) + g(\widehat X_2)
- 1
\end{equation}

By Definition \ref{dfn:hat} $X_2$ is obtained from $\widehat X_2$
by Dehn filling, and moreover, the core of the attached solid
torus is parallel into a meridian curve on $\del X_2$. It is
therefore parallel into every Heegaard surface of $X_2$. In
\cite{rieck} this type of filling is called \em good \em and it
was shown there that for good fillings:

\begin{equation}\label{equ:dehn-surgery}
g(\widehat X_2) - 1 \leq g(X_2) \leq g(\widehat X_2)
\end{equation}

Equations (\ref{equ:amalgamation}) and (\ref{equ:dehn-surgery})
complete the proof of Theorem \ref{thm:bound-on-degeneration}.
\end{proof}

\section{Morimoto's Conjecture for $n=2$}
\label{sec:morimoto}

\begin{thm}
\label{thm:morimoto-n=2}

Suppose $X_1$ and $X_2$ are m-small knots and that $g(X)$ $<
g(X_1) + g(X_2)$.  Then at least one of $X_1$ or $X_2$ contains a
\pmm.
\end{thm}

\begin{rmkks}
\label{rmk:case-n=2} {\rm
\begin{enumerate}
\item Morimoto proved Theorem \ref{thm:morimoto-n=2} for m-small
knots in $S^3$ (Theorem 1.6 of \cite{morimoto1}).

\item We again emphasize that there is no symmetry between $X_1$
and $X_2$. The the swallow follow torus found in Theorem
\ref{thm:basic} follows one (say $K_1$) and swallows the other
($K_2$), and as we shall see in the proof, in this case $X_2$
contains a \pmm. By Moriah and Rubinstein \cite{moriah-rubinstein}
and Morimoto, Sakuma and Yokota \cite{skuma-morimoto-yokata} there
exist knots without a \pmm. Taking such knot as $K_1$ and
connecting it with $K_2$ a knot that does have a \pmm\ (for
example any 2-bridge knot, by Kobayashi \cite{kobayashi-2-bridge})
we always get $X = X_1 \cup_T \widehat X_2$ and never $X =
\widehat X_1 \cup_T X_2$. This is used in Section
\ref{sec:induction}.

\item Note that the proof in this section holds in greater
generality than stated: all we need is that $X_2$ is the exterior
of a prime knot, but $X_1$ may be the exterior of a connected sum
of many components. This too is used in Section
\ref{sec:induction} (as the induction hypothesis).
\end{enumerate}
}
\end{rmkks}

\begin{proof}[Sketch of proof of Theorem \ref{thm:morimoto-n=2}]

Let $X=X_1 \star X_2$ and let $\s$ be a minimal genus Heegaard
splitting for $X$. By Theorem \ref{thm:basic} $\s$ weakly reduces
to a \sft\ (say $T$). Thus $X = X_1 \cup_T \widehat X_2$ or $X =
\widehat  X_1 \cup_T X_2$, say the former.  Our goal is showing
that $X_2$ contains a \pmm. This follows from the following three
statements.  In all three we assume that $X$ is the exterior of a
connected sum of m-small knots; the assumption that $X_2$ is prime
is used only in Theorem \ref{thm:intersecting-annulus-once} (we do
not assume that $X_1$ is prime).

\begin{lem}
\label{lem:widehat-x2=x2}

If $g(X) < g(X_1) + g(X_2)$ then $g(X_2) = g(\widehat X_2).$
\end{lem}

\begin{thm}
\label{thm:intersecting-annulus-once}

Let $A$ be the essential annulus in $\widehat X_2$ with one
boundary component a meridian of $X_2$ and the other a longitude
of $\del \widehat X_2 \setminus \del X_2$.

Then there is a minimal genus \hhs\ for $\widehat X_2$ that
intersects $A$ in exactly one essential \scc.
\end{thm}

\begin{rmk}
Note that Theorem \ref{thm:intersecting-annulus-once} implies that
$\widehat X_2$ has a minimal genus Heegaard surface that separates
the boundary components.  We do not know if this is always the
case.
\end{rmk}

By Definition \ref{dfn:hat} $X_2$ is obtained from $\widehat X_2$
by Dehn filling; therefore any \hhs\ for $\widehat X_2$ is a \hhs\
for $X_2$.

\begin{lem}
\label{lem:there-is-pm}

The Heegaard surface found in Theorem
\ref{thm:intersecting-annulus-once} (when considered as a \hhs\
for $X_2$) contains a \pmm.
\end{lem}

Theorem \ref{thm:morimoto-n=2} clearly follows since by Lemma
\ref{lem:widehat-x2=x2} the surface found in Theorem
\ref{thm:intersecting-annulus-once} is minimal genus in $X_2$. We
conclude this section by indicating the proofs of these statements
in the order in which they appeared:

\begin{proof}[Proof of Lemma \ref{lem:widehat-x2=x2}]

From the previous section, recall Equations
(\ref{equ:amalgamation}): $g(X) = g(X_1) + g(\widehat X_2) - 1$,
and (\ref{equ:dehn-surgery}): $g(\widehat X_2) - 1 \leq g(X_2)
\leq g(\widehat X_2).$ The lemma follows.
\end{proof}

\begin{proof}[Sketch of proof for Theorem \ref{thm:intersecting-annulus-once}]

This is the most difficult part of the proof and accordingly our
sketch is very rough.  Let $\widehat \s_2 \subset \widehat X_2$ be
a minimal genus Heegaard surface.  The proof has three cases not
unlike the proof of Theorem \ref{thm:basic}.

\begin{description}

\item[Case one: assume $\widehat \s_2$ is strongly irreducible]
Since $A$ is essential and we assumed $\s$ to be strongly
irreducible we can isotope $\s$ to intersect $A$ is essential
curves only. Denote the subannuli of $A$ cut open along $\s$ by
$\{A_j\}_{j=0}^k$, with $A_0$ adjacent to $\del X_2$ and $A_k$
adjacent to $\del \widehat X_2 \setminus \del X_2$. Next we fill
$\del \widehat X_2 \setminus \del X_2$ to obtain $X_2$, and denote
the core of the attached solid torus by $\gamma$.  Note that we
can still consider $A$ as an embedded annulus (although not
properly) with one curve of $\del A$ a meridian of $X_2$ and the
other $\gamma$. We show that the equator of $A_1$ is a core of a
1-handle (of one of the compression bodies obtained by cutting
$X_2$ open along $\widehat \s_2$) by compressing $\widehat \s_2$
maximally into the side containing $A_1$ but only along disks
disjoint from $A_1$. By applying Scharlemann's No Nesting Lemma
\cite{no-nestig}, we see that either the surface we obtain is a
torus bounding a solid torus with the core curve of $A_1$ as its
core, or the components of this surface cut open along $A$ are
essential.  Since by cutting $\widehat X_2$ along $A$ we obtain a
manifold homeomorphic to $X_2$, the latter contradicts
m-smallness.  The former proves the theorem by isotoping $\gamma$
along $A$ to the core of $A_1$ and drilling it out again.

\end{description}

From now on we assume $\widehat \s_2$ weakly reduces. We
untelescope (see Section \ref{sec:background}) $\widehat \s_2$ to
obtain a collection of essential surfaces $S_i$ and consider two
more cases:

\begin{description}

\item[Case Two: $(\forall j)$ $S_j \cap A = \emptyset$] let $M$ be
the component of $\widehat X_2$ cut open along $S=\cup_i S_i$
containing $A$. Since the induced Heegaard splitting is strongly
irreducible \cite{scharl-abby} as in Case One we can arrange for
it to intersect $A$ in a single essential curve.  By amalgamating,
we retrieve $\widehat \s_2$.  We complete the proof by using an
explicit description of amalgamation to verify that this can be
done without introducing curves of intersection with $A$.

\item[Case Three: $(\exists j)$ $F_j \cap A \neq \emptyset$] (this
is similar to an argument of Morimoto \cite{morimoto2}) each
component of $F_j$ cut open along $A$ is incompressible.  It is
therefore either essential or a boundary compressible annulus
(since for surfaces of negative Euler characteristic boundary
compression implies compression).  Thus $F_j$ is constructed by
pasting annuli together, and is a torus.  By cutting $\widehat
X_2$ and $F_j$ along $A$ we obtain meridional incompressible
annuli in $X_2$, that must be boundary parallel by our assumption
that $X_2$ is a prime knot.  This implies that $F_j$ is a boundary
parallel torus in $\widehat X_2$.  However, a surface obtained by
untelescoping a minimal genus \hhs\ is never boundary parallel.
With contradiction we see that Case Three never occurs.
\end{description}

This completes our sketch the proof of Theorem
\ref{thm:intersecting-annulus-once}.
\end{proof}

\begin{proof}[Proof of Lemma \ref{lem:there-is-pm}]

This is immediate.  After filling $\widehat X_2$ to get $X_2$ the
meridian of the attached solid torus extends to a meridian of the
handlebody that intersects the meridional annulus $A$ exactly
once.
\end{proof}

This completes our sketch of the proof of Theorem
\ref{thm:morimoto-n=2}.
\end{proof}

\section{Proof of Morimoto's Conjecture for m-small knots}
\label{sec:induction}

The Swallow Follow Torus Theorem (\ref{thm:basic}) and the
Numerical Theorem (\ref{thm:numerical}) were proved for the
connected sum of any number of m-small knots.  We conclude this
paper by the inductive step of the proof of Morimoto's Conjecture:

\begin{qthm}[\ref{thm:geometric}]
Let $\{X_i\}_{i=1}^n$ be exteriors of m-small knots. Then:

If $g(X) < \s_{i=1}^n g(X_i)$ there exists some $I \subset
\{1,..,n\}$ a proper subset so that $X = \star_{i \in I} X_i$
contains a \pmm.
\end{qthm}

\begin{proof}[Sketch of proof.]

Let $\s$ be a minimal genus Heegaard surface for $X$. By Theorem
\ref{thm:basic} $\s$ weakly reduces to a \sft\ $T$ yielding the
decomposition $X_I \cup_T \widehat X_J$ (with $I \sqcup J
=\{1,..,n\}, I,J \neq \emptyset$).  We assign a complexity to a
knot exterior $Y$: the number of prime summands, denoted $n(Y)$.
To the decomposition $X = X_I \cup_T \widehat X_J$ we assign the
complexity $(n(X),n(X_J))$ (ordered lexicographically). Suppose $X
= X_I \cup_T \widehat X_J$ is a minimal complexity counterexample
to Theorem \ref{thm:geometric}.

If $g(X) = g(X_I) + g(X_J)$ then, since $g(X) < \s_{i=1}^n
g(X_i)$, either $g(X_I) < \s_{i \in I} g(X_i)$ or $g(X_J) < \s_{i
\in J} g(X_i)$.  By the minimality assumption, in the former case
some subcollection of $I$ has a \pmm\ and in the latter some
subcollection of $J$.  In either case $X$ is not a counterexample.
This contradiction together with Theorem
\ref{thm:bound-on-degeneration} shows:

   $$g(X) = g(X_I) + g(X_J) - 1.$$

If $n(X_J)=1$ by (3) of Remark \ref{rmk:case-n=2}, $X_J$ contains
a \pmm. So we may assume that $n(X_J)>1$. By Theorem
\ref{thm:basic} (applied to $X_J$) the minimal genus Heegaard
splitting of $X_J$ induced by the minimal genus Heegaard splitting
of $X$ weakly reduces to a \sft\ in $X_J$, say $X_J = X_{J_1}
\cup_{T_J} \widehat X_{J_2}$. If $g(X_J) < g(X_{J_1}) +
g(X_{J_2})$ then (by the minimality of $(n(X),n(X_J))$) some
subcollection of $X_J$ has a \pmm\ and $X$ is not a
counterexample, a contradiction. We get:

     $$g(X_J) = g(X_{J_1}) + g(X_{J_2}).$$

Next we consider $X_I \cup_T X_{J_1}$.  If $g(X_I \cup_T X_{J_1})
< g(X_I) + g(X_{J_1})$ then some subcollection of $\{X_i\}_{i \in
I \cup J_1}$ contains a \pmm\ (minimality again)  so $X$ is not a
counterexample, a contradiction. We get:

 $$g(X_I \cup_T X_{J_1}) = g(X_I) + g(X_{J_1}).$$

Since the Heegaard splitting induced by $\s$ on $X_J$ weakly
reduces to $T_J$, $\s$ also weakly reduces to $T_J$.  This gives
$(X_I \cup_T \widehat X_{J_1}) \cup_{T_J} \widehat X_{J_2}$.  The
equations above imply that $g(X) = g(X_I \cup_T \widehat X_{J_1})
+ g(X_{J_2}) - 1$.  By the minimality, $(X,X_{J_2})$ is not a
counterexample, and we conclude that $X$ itself is not a
counterexample at all. With this contradiction we complete the
proof.
\end{proof}

% ----------------------------------------------------------------

\providecommand{\bysame}{\leavevmode\hbox
to3em{\hrulefill}\thinspace}
\providecommand{\MR}{\relax\ifhmode\unskip\space\fi MR }
\providecommand{\MRhref}[2]{%
  \href{http://www.ams.org/mathscinet-getitem?mr=#1}{#2}
} \providecommand{\href}[2]{#2}

\end{document}